\documentclass[runningheads]{llncs}
\usepackage{url}
\usepackage{multirow}
\usepackage{amssymb}
\usepackage{amsxtra}
\usepackage{makeidx}
\usepackage{amsfonts}
\usepackage{algorithm}
\usepackage{algorithmic}
\usepackage{float}

\begin{document}
\title{Mirror Descent and Constrained Online Optimization Problems\thanks{The research by Alexander A. Titov and Fedor S. Stonyakin (Theorem 2 and Remark 3) was partially supported by Russian Science Foundation according to the research project 18-71-00048.}}
%
%
\author{Alexander A. Titov\inst{1}\and
Fedor S. Stonyakin\inst{1,2} \and
Alexander V. Gasnikov\inst{1} \and Mohammad S. Alkousa\inst{1}}
\authorrunning{A.A.~Titov et al.}
%
\institute{Moscow Institute of Physics and Technologies, Moscow\\
	\email{a.a.titov@phystech.edu, gasnikov@yandex.ru, mohammad.alkousa@phystech.edu},\\
	\and
	V.\,I.\,Vernadsky Crimean Federal University, Simferopol\\
	\email{fedyor@mail.ru}}
\maketitle
\begin{abstract}
We consider the following class of online optimization problems with functional constraints. Assume, that a finite set of convex Lipschitz-continuous non-smooth functionals are given on a closed set of $n$-dimensional vector space. The problem is to minimize the arithmetic mean of functionals with a convex Lipschitz-continuous non-smooth constraint. In addition, it is allowed to calculate the (sub)gradient of each functional only once. Using some recently proposed adaptive methods of Mirror Descent the method is suggested to solve the mentioned constrained online optimization problem with optimal estimate of accuracy. For the corresponding non-Euclidean prox-structure the case of a set of $n$-dimensional vectors lying on the standard $n$-dimensional simplex is considered.
\end{abstract}

\section{Introduction}
Online convex optimization plays a key role in solving the problems, where statistical information is being updated \cite{Hazan-Kale,Hazan-2016}. There are a lot of examples of such problems, concerning internet network, consumer data sets or financial market. Quite a few branches of science also face the above mentioned problems, for example machine learning applications  \cite{Jenatton_adaptive_with_constraints}. The important example is the descision-making problem  \cite{Hazan-2016,Kalai_efficient_algorithms}. Suppose, we are given $N$ experts and range of admissible solutions lie on the unit simplex. Every expert gives his estimates of losses with the possible solution and the problem is to minimize total losses from the point view of all experts (the arithmetic mean). Therefore, in recent years, methods for solving online optimization problems have been actively developed \cite{Bubeck-Eldan,Bubeck-Bianchi,Gasnikov2015,Gasnikov2017,Hazan-Kale,Hazan-2016,Jenatton_adaptive_with_constraints,Lugosi-Bianchi}.

In problems of online convex optimization, it is required to minimize the sum (or the arithmetic mean) of several convex Lipschitz functionals $f_i$ ($i = \overline{1, N}$) given on some closed set $ Q \subset \mathbb{R}^n $. It should be noted that it is possible to calculate the (sub)gradient $\nabla f_i(x)$ of each functional $f_i$ only once. Our paper is devoted to some optimal methods for the following type of problems
\begin{equation}\label{Problem_Statement}
\begin{cases}
\frac{1}{N}\sum\limits_{i=1}^N f_i(x) \rightarrow \min\limits_{x \in Q}\\
s.t. \ \  g(x)\leq 0
\end{cases}
\end{equation}

We assume that the functionals $f_i$ and $g$ satisfy the Lipschitz property, i.e. there exists a number $M > 0$, such that
\begin{equation}
|g(x)-g(y)|\leq M\|x-y\|,
\end{equation}
\begin{equation}
|f_i(x)-f_i(y)|\leq M\|x-y\| \quad \forall i = \overline{{1, N}}.
\end{equation}

We can explain the meaning of such formulation of the problem in the following situation. Suppose that we are engaged in some kind of activity during the fixed number of days. Each day can be productive or non-productive. We want to live out $N$ productive days (not necessarily in a row, there can be some non-productive days within this period), so that the total nerve costs (characterized by $f_i (x)$) would be minimal. Note that we pay nervous expenses only in productive days, when we try to do something. In non-productive days we do nothing, our aim is to return to the productive state, but we do not pay any costs. The productivity of the day is determined by the condition $g(x^k) \le \varepsilon$. Let's define index $i$ as the number of the productive day. This day we receive feedback from the outside world in the next form: $\nabla f_i(x^k)$ and using this information we build a strategy for the next day $x^{k + 1}$. In non-productive days, we get information about how far have we gone out of the functional constraint and we try to return to this framework. There is no point in arranging unnecessary non-productive days. Therefore, it is also desirable to minimize the number of non-productive days for a given $N$. The proposed algorithm provides a small amount of costs simultaneously, ensuring that the number of non-productive days will be no more than $O(N)$.

The optimization problems of non-smooth functionals with constraints attract widespread interest in large-scale optimization and its applications \cite{bib_ttd,bib_Shpirko}. There are various methods of solving this kind of optimization problems. Some examples of these methods are: bundle-level method \cite{bib_Nesterov}, penalty method \cite{bib_Vasilyev}, Lagrange multipliers method \cite{bib_Boyd}. Among them, Mirror Descent (MD) \cite{beck2003mirror,nemirovsky1983problem} is viewed as a simple method for non-smooth convex optimization.

Note that a functional constraint, generally, can be non-smooth. That is why we consider subgradient methods. These methods have a long history starting from the method for deterministic unconstrained problems and Euclidean setting in \cite{shor1967generalized} and the generalization for constrained problems in \cite{polyak1967general}, where the idea of steps switching between the direction of subgradient of the objective and the direction of subgradient of the constraint was suggested. Non-Euclidean extension, usually referred to as Mirror Descent, originated in \cite{nemirovskii1979efficient,nemirovsky1983problem} and was later analyzed in \cite{beck2003mirror}. An extension for constrained problems was proposed in \cite{nemirovsky1983problem}, see also a recent version in \cite{beck2010comirror}.

Usually, the stepsize and stopping rule for Mirror Descent requires to know the Lipschitz constant of the objective function and constraint, if any. Adaptive stepsizes, which do not require this information, are considered in \cite{bib_Nemirovski} unconstrained problems, and in \cite{beck2010comirror} for constrained problems. Recently, in \cite{bib_Adaptive} optimal algorithms of Mirror Descent for convex programming problems with Lipschitz functional constraints with both adaptive step selection and adaptive stopping criteria were proposed for a number of classes of problems. Also there were considered some modifications of these methods for the case of problems with many functional constraints in \cite{bib_Stonyakin}. In \cite{Jenatton_adaptive_with_constraints} authors considered adaptive algorithms for online convex optimization problem with Constraints, but with only standard Euclidean prox-structure.

In this paper we propose adaptive and non-adaptive algorithms for solving the problem \eqref{Problem_Statement}. Note that we consider arbitrary proximal structure, which seems essential for the problem of experts \cite{Gasnikov2015,Gasnikov2017,Hazan-Kale,Hazan-2016}. The paper consists of Introduction and five main sections. In Section 2 we give some basic notation concerning convex optimization problems with functional constrains and online optimization problems. In section 3 we propose a non-adaptive algorithm of Mirror Descent for the considered online optimization problem \eqref{Problem_Statement}. Section 4 is devoted to an adaptive analog of this method (Algorithm \ref{alg02}).

Also in section 4, by analogy with \cite{bib_Stonyakin}, we propose a modification of Algorithm \ref{alg02} for problems with several functional constraints (Algorithm \ref{alg03}). It is shown that Algorithms \ref{alg01}, \ref{alg02} and \ref{alg03} are optimal accurate to multiplication by constants under the condition of nonnegativity of the regret (see Theorems 1 and 2). In section 5 the condition of negative regret is considered. In this case we get the optimal quality of estimation by the objective function, but the estimation of the number of non-productive steps is worse than \eqref{Non_neg_nonprod}. In the last section we consider some numerical experiments that allow us to compare the work of Algorithms \ref{alg01}, \ref{alg02}, and \ref{alg03} for certain examples.

Summing up, contributions of this paper are as follows:
\begin{itemize}
	\item two methods (adaptive and non-adaptive) were proposed to solve the online optimization problem for an arbitrary prox-structure;
	\item the number of non-productive steps is $O(N)$ in the case of nonnegative regret;
	\item the number of non-productive steps is $O(N^2)$, but the accuracy by regret is better.
\end{itemize}

\section{Problem Statement and Standard Mirror Descent Basics}

Let $(E,||\cdot||)$ be a normed finite-dimensional vector space and $E^*$ be the conjugate space of $E$ with the norm:
$$||y||_*=\max\limits_x\{\langle y,x\rangle,||x||\leq1\},$$
where $\langle y,x\rangle$ is the value of the continuous linear functional $y$ at $x \in E$.

Let $Q\subset E$ be a (simple) closed convex set, $d : Q \rightarrow \mathbb{R}$ be a distance generating function (d.g.f) which is continuously differentiable and $1$-strongly convex w.r.t. the norm $\|\cdot\|$, i.e.
$$\forall x, y \in Q \hspace{0.2cm} \langle \nabla d(x) - \nabla d(y), x-y \rangle \geq \| x-y \|^2,$$
and assume that $\min\limits_{x\in Q} d(x) = d(0).$
Suppose, we have a constant $\Theta_0$ such that $d(x_{*}) \leq \Theta_0^2,$ where $x_*$ is a solution of \eqref{Problem_Statement}.

Note that if there is a set of optimal points for \eqref{Problem_Statement} $X_* \subset Q$,  we may assume that
$$\min\limits_{x_* \in X_*} d(x_*) \leq \Theta_0^2.$$
For all $x, y\in Q \subset E$ consider the corresponding Bregman divergence
$$V(x, y) = d(y) - d(x) - \langle \nabla d(x), y-x \rangle.$$
Standard proximal setups, i.e. Euclidean, entropy, $\ell_1/\ell_2$, simplex, nuclear norm, spectahedron can be found, e.g. in \cite{bib_Nemirovski}. Let us define the proximal mapping operator standardly
$$
\mathrm{Mirr}_x (p) = \arg\min\limits_{u\in Q} \big\{ \langle p, u \rangle + V(x, u) \big\} \;  \text{ for each } \;  x\in Q \; \text{ and } \;  p\in E^*.
$$
We make the simplicity assumption, which means that $\mathrm{Mirr}_x (p)$ is easily computable. There are well-known examples of distance generating function, let us denote $\ell_p$ norm by $\|x\|_{p} $, and the unit simplex in $\mathbb{R}^n$ by
$$S_n(1) = \left\{ x\in \mathbb{R}_{+}^n \;|\; \displaystyle \sum_{i=1}^{n}x_i = 1\right\}.$$ Consider two cases:
\begin{itemize}
	\item  if  $p=1$, then
	\begin{equation}\label{neevklid_prox_1}
	d(x) = \ln n + \displaystyle \sum_{k=1}^{n}x_k \ln x_k, \quad
	V(x,y) = \displaystyle \sum_{k=1}^{n}x_k \ln\left(\frac{x_k}{y_k}\right);
	\end{equation}
	\item  if $p=2$, then $d(x) = \frac{1}{2} \|x\|_{2}^{2} ,$  $V(x,y) = \frac{1}{2} \|x-y\|_{2}^{2}.$
\end{itemize}

Let $ Q = \mathbb{B}_p^{n}(1) = \{ x \in \mathbb{R}^{n}; \|x\|_{p} \leq 1 \}$ be the unit ball with $l_p$ norm. One can note the following: if $p \geq 2 $, then it is optimal to choose the $l_2$-norm and the Euclidean prox-structure.

Define $q$ by $\frac{1}{p} + \frac{1}{q} = 1$ and consider $ 1 \leq p \leq 2$, then $q \geq 2$. If in this case $q = O(\ln n)$, then it is optimal to choose $l_p$-norm and prox-structure with distance generating function $$d(x) = \frac{1}{2(p-1)}\|x\|_{p}^{2}.$$ In all these cases $R^2 = \max\limits_{x\in Q} d(x) \geqslant \Theta_0^2$.

For $q > \Omega (\ln n)$, we choose $l_a$-norm, where $$a= \frac{2 \ln n}{2 \ln n -1}$$ and prox-structure with distance generating function $$d(x) = \frac{1}{2(a-1)}\|x\|_{a}^{2}.$$ In this case
\begin{equation}\label{neevklid_prox_2}
R^2 = O(\ln n) \geqslant \Theta_0^2 \; \text{ and } \; \Theta_0 \leqslant O(\sqrt{\ln n}).
\end{equation}

Let us remind one well-known statement (see, e.g. \cite{bib_Nemirovski}).

\begin{lemma}\label{lem1}
	Let $f: Q \rightarrow \mathbb{R}$ be a convex subdifferentiable function over the convex set $Q$ and  $z=Mirr_{y}(h \nabla f(y))$ for some $h>0$, $y, z \in Q$. Then for each $x\in Q$
	\begin{equation}\label{eq7}
	h\langle\nabla f(y), y-x\rangle\leq\frac{h^2}{2}||\nabla f(y)||_*^2 + V(y,x) - V(z,x).
	\end{equation}
\end{lemma}

\section{Online Optimization for the Case of Non-negative Regret: Non-Adaptive Algorithm}

Assume that the method produces $N$ productive steps and each step the (sub)gradient of exactly one functional of the objectives is calculated. Denote the number of non-productive steps by $N_J$. Let's consider the non-adaptive method for the problem \eqref{Problem_Statement} with a constant step, which depends on the Lipschitz constant $M$. As a result, we get a sequence
$\{x^k \}_{k \in I}$ (on productive steps), which can be considered as a solution to the problem \eqref{Problem_Statement}  with accuracy $ \delta $ (see \eqref{eqiv_crit}).

\begin{algorithm}
	\caption{Constrained Online Optimization: Non-Adaptive Mirror Descent Algorithm}
	\label{alg01}
	\begin{algorithmic}[1]
		\REQUIRE $\varepsilon, N, \Theta_0^2,Q,d(\cdot), x^0$
		\STATE $i:=1, \; k:=0;$
		\REPEAT
		\IF{$g(x^k) \leqslant \varepsilon$}
		\STATE $h = \frac{\varepsilon}{M^2};$
		\STATE $x^{k+1}:=Mirr[x^k](h\nabla f_i(x^k));$
		\STATE $i:=i+1;$
		\STATE $k:=k+1;$
		\ELSE
		\STATE $h= \frac{\varepsilon}{M^2};$
		\STATE $x^{k+1}:=Mirr[x^k](h\nabla g(x^k));$
		\STATE $k:=k+1;$
		\ENDIF
		\UNTIL{$i=N+1$}
		\STATE Guaranteed accuracy:
		\begin{equation}\label{eqiv_crit}
		\delta:= \frac{\varepsilon}{2} + \frac{M^2\Theta_0^2 }{\varepsilon N} -\frac{\varepsilon N_J}{2N}
		\end{equation}
	\end{algorithmic}
\end{algorithm}

By Lemma \ref{lem1}
$$f_i(x^k)-f_i(x)\leq \frac{h}{2}M^2+\frac{V(x^k,x)}{h}-\frac{V(x^{k+1},x)}{h} = \frac{\varepsilon}{2} + \frac{V(x^k,x)}{h}-\frac{V(x^{k+1},x)}{h}$$
$$g(x^k)-g(x)\leq \frac{h}{2}M^2+\frac{V(x^k,x)}{h}-\frac{V(x^{k+1},x)}{h} = \frac{\varepsilon}{2} + \frac{V(x^k,x)}{h}-\frac{V(x^{k+1},x)}{h} $$
Taking summation over productive and non-productive steps, we get
$$
\sum\limits_{i = 1}^{N}(f_i(x^k)-f_i(x^*))+\sum\limits_{k\in J}(g(x^k)-g(x^*))\leq
$$
$$\frac{\varepsilon}{2}(N + N_J)+ \frac{1}{h}\sum\limits_{k = 0}^{N + N_J - 1}\left(V(x^k,x^*)-V(x^{k+1},x^*)\right) =
$$
$$
= \frac{\varepsilon}{2}(N + N_J)+ \frac{M^2}{\varepsilon}\sum\limits_{k = 0}^{N + N_J - 1}\left(V(x^k,x^*)-V(x^{k+1},x^*)\right),
$$
then
\begin{equation}\label{eq5}
\sum\limits_{i=1}^{N}(f_i(x^k)-f_i(x^*)) \leq \frac{\varepsilon}{2}N + \frac{M^2\Theta_0^2 }{\varepsilon} -\frac{\varepsilon}{2}N_J
\end{equation}
and by virtue of \eqref{eqiv_crit}
\begin{equation}\label{eq6}
\frac{1}{N}\sum\limits_{i=1}^Nf_i(x^k)-\min\limits_{x\in Q}\frac{1}{N}\sum\limits_{i=1}^Nf_i(x) \leq \delta.
\end{equation}
If we assume the nonnegativity of the regret (i.e. the left side in \eqref{eq5}) and
\begin{equation}\label{Accur}
\delta \leq \varepsilon = \frac{C}{\sqrt{N}} \text{ for some } C > 0,
\end{equation}
then we get
$$
0 \leq N + \frac{2M^2\Theta_0^2}{\varepsilon^2} - N_J = N + \frac{2M^2\Theta_0^2}{C^2}N - N_J,
$$
then
$$
N_J \leq N \cdot\left(1 + \frac{2M^2\Theta_0^2}{C^2}\right) \sim O(N).
$$

Thus, we have the following result
\begin{theorem}
	Suppose Algorithm 1 works exactly $N$ productive steps. After the stopping of the Algorithm 1, the following inequality holds:
	$$
	\frac{1}{N}\sum\limits_{i=1}^Nf_i(x^k)-\min\limits_{x\in Q}\frac{1}{N}\sum\limits_{i=1}^N f_i(x) \leq \delta.
	$$
	For the case \eqref{Accur} and
	$$
	\frac{1}{N}\sum\limits_{i=1}^Nf_i(x^k)-\min\limits_{x\in Q}\frac{1}{N}\sum\limits_{i=1}^Nf_i(x) \geq 0
	$$
	there will be no more than
	\begin{equation}\label{estim_nonprod}
	N \cdot\left(1 + \frac{2M^2\Theta_0^2}{C^2}\right) \sim O(N).
	\end{equation}
	non-productive steps.
\end{theorem}

\begin{remark}
	The estimate \eqref{estim_nonprod} is optimal for the considered class of problems \cite{Hazan-Kale}.
\end{remark}

\begin{corollary}
	If $Q = S_n(1)$ and the corresponding prox-structure is chosen as \eqref{neevklid_prox_1}, then by \eqref{neevklid_prox_2} the estimate \eqref{estim_nonprod} modifies into
	$$
	N_J \leqslant N \cdot \left(1 + \frac{2M^2 \ln n}{C^2}\right).
	$$
\end{corollary}

\section{Adaptive Mirror Descent for the Case of Non-negative Regret}

Now, let us consider the adaptive analog of Algorithm 1 for problem \eqref{Problem_Statement}. The main feature is a nondecreasing stepsize with consideration of the norm of (sub)gradient of the objective function or the constraints in a particular step. Therefore, the proposed algorithm will work until there are exactly $N$ productive steps. As a result, we get a sequence $\{x^k \}_{k \in I}$ on productive steps, which can be considered as a solution to the problem \eqref{Problem_Statement} with accuracy $ \delta $ (see \eqref{eqiv_crit_1}).

\begin{algorithm}[t]
	\caption{Constrained Online Optimization: Adaptive Mirror Descent Algorithm}
	\label{alg02}
	\begin{algorithmic}[1]
		\REQUIRE  $\varepsilon, N, \Theta_0^2,Q,d(\cdot), x^0: \text{ and } \sup\limits_{x, y \in Q} V(x, y) \leqslant \Theta_0^2$ 
		\STATE $i:=1, \; k:=0;$
		\REPEAT
		\IF{$g(x^k) \leqslant \varepsilon$}
		\STATE $M_k := \|\nabla f_i(x^k)\|_*;$
		\STATE $h_k = \Theta_0 \left(\sum\limits_{t = 0}^k M_t^2\right)^{-1/2};$
		\STATE $x^{k+1}:=Mirr[x^k](h_k\nabla f_i(x^k));$
		\STATE $i:=i+1;$
		\STATE $k:=k+1;$
		\ELSE
		\STATE $M_k := \|\nabla g(x^k)\|_*;$
		\STATE $h_k = \Theta_0 \left(\sum\limits_{t = 0}^k M_t^2\right)^{-1/2};$
		\STATE $x^{k+1}:=Mirr[x^k](h_k\nabla g(x^k));$
		\STATE $k:=k+1;$
		\ENDIF
		\UNTIL{$i=N+1$}
		\STATE Guaranteed accuracy:
		\begin{equation} \label{eqiv_crit_1}
		\delta:= \frac{2\Theta_0}{N} \left(\sum\limits_{i=0}^{N+N_J-1} M_i^2\right)^{1/2} - \varepsilon \cdot \frac{N_J}{N}.
		\end{equation}
	\end{algorithmic}
\end{algorithm}

By Lemma \ref{lem1}
$$f_i(x^k)-f_i(x)\leq \frac{h_k}{2}\|\nabla f_i(x^k)\|_*^2+\frac{V(x^k,x)}{h_k}-\frac{V(x^{k+1},x)}{h_k}$$
$$g(x^k)-g(x)\leq \frac{h_k}{2}\|\nabla g(x^k)\|_*^2+\frac{V(x^k,x)}{h_k}-\frac{V(x^{k+1},x)}{h_k} $$

Dividing each inequality by $h_k$ and summing up for $k$ from $0$ to $N+N_J-1$, and by using the definition of $h_k$, we obtain
$$
\sum\limits_{k\in I}  \big( f(x^k) - f(x_{*}) \big) + \sum\limits_{k\in J}  \big( g(x^k) - g(x_{*}) \big) \leq \sum\limits_{k = 0}^{N+N_J-1} \frac{h_k M_k^2}{2} +
$$
$$+ \sum\limits_{k = 0}^{N+N_J-1}\frac{1}{h_k} \left(V(x^k,x_{*}) - V(x^{k+1},x_{*}) \right) \text{ and }
$$
$$	
\sum\limits_{k = 0}^{N+N_J-1} \frac{1}{h_k} \big(V(x^k,x_{*}) - V(x^{k+1}, x_{*})  \big)
= \frac{1}{h_0} V(x^0, x_*) + \sum\limits_{k=0}^{N+N_J-2} \Big(\frac{1}{h_{k+1}} - \frac{1}{h_k} \Big) V(x^{k+1}, x_*) -
$$
$$
- \frac{1}{h_{N+N_J-1}} V\left(x^{N+N_J}, x_*\right) \leq \frac{\Theta_0^2}{h_0} + \Theta_0^2 \sum\limits_{k=0}^{N+N_J-2} \Big(\frac{1}{h_{k+1}} - \frac{1}{h_k} \Big) = \frac{\Theta_0^2}{h_{N+N_J-1}}.
$$
Whence, by the definition of stepsizes $h_k$,
\begin{align}
&\sum\limits_{i = 1}^{N}  \big(f_i(x^k) - f(x_{*}) \big) + \sum\limits_{k\in J}  \big(g(x^k) - g(x_{*}) \big) \leq  \sum\limits_{k =  0}^{N+N_J-1} \frac{h_k M_k^2}{2} +\frac{\Theta_0^2}{h_{N+N_J-1}} \notag \\
& \leq \sum\limits_{k=0}^{N+N_J-1} \frac{\Theta_0}{2} \frac{M_k^2}{\left( \sum_{j=0}^k M_j^2 \right)^{1/2}} + \Theta_0\left( \sum_{k=0}^{N+N_J-1} M_k^2 \right)^{1/2} \leq 2 \Theta_0 \left(\sum_{k=0}^{N+N_J-1} M_k^2 \right)^{1/2}
\end{align}
where we used inequality $$\sum\limits_{i=0}^{N+N_J-1}  \frac{M_i^2}{\left(\sum_{j=0}^i M_j^2 \right)^{1/2}} \leq 2\left( \sum_{i=0}^{N+N_J-1} M_i^2 \right)^{1/2},$$
which can be proved by induction. Since, for $k\in J$, $g(x^k) - g(x_{*}) \geq g(x^k) > \varepsilon$, we get
\begin{equation}\label{Adaptive_estimate}
\sum\limits_{i=1}^{N}(f_i(x^k)-f_i(x^*))  < \varepsilon N - \varepsilon(N +N_J) + 2\Theta_0 \left(\sum\limits_{i=0}^{N+N_J-1} M_i^2 \right)^{1/2}.
\end{equation}
and by \eqref{eqiv_crit_1}
\begin{equation}\label{eq6}
\frac{1}{N}\sum\limits_{i=1}^Nf_i(x^k)-\min\limits_{x\in Q}\frac{1}{N}\sum\limits_{i=1}^Nf_i(x) \leq \delta.
\end{equation}
If we assume the nonnegativity of the regret (i.e. the left side in \eqref{Adaptive_estimate}) and the accuracy is given by \eqref{Accur},
one can get
$$
\varepsilon(N +N_J) \leq \varepsilon N + 2\Theta_0 \left(\sum\limits_{i=0}^{N+N_J-1} M_i^2 \right)^{1/2} \leq \varepsilon N + 2M\Theta_0\cdot\sqrt{N+N_J},
$$
$$
N_J^2 \leq \frac{4M^2\Theta_0^2 (N+N_J)}{\varepsilon^2} = \frac{4M^2\Theta_0^2 (N+N_J)N}{C^2}
$$
Further,
$$
\frac{N_J^2}{N^2 + NN_J} = \frac{\left(\frac{N_J}{N}\right)^2}{1 + \frac{N_J}{N}} \leq \frac{4M^2\Theta_0^2}{C^2}
$$
and $N_J = O(N)$. Thus, we have come to the following result.

\begin{theorem}
	Suppose Algorithm 2 works exactly $N$ productive steps. After the stopping of the Algorithm 2, the following inequality holds:
	$$
	\frac{1}{N}\sum\limits_{i=1}^Nf_i(x^k)-\min\limits_{x\in Q}\frac{1}{N}\sum\limits_{i=1}^Nf_i(x) \leq \delta.
	$$
	For the case of \eqref{Accur} and
	$$
	\frac{1}{N}\sum\limits_{i=1}^Nf_i(x^k)-\min\limits_{x\in Q}\frac{1}{N} \sum\limits_{i=1}^Nf_i(x) \geq 0
	$$
	there will be no more than $O(N)$ non-productive steps.
\end{theorem}

\begin{remark}
	Algorithm 2 is optimal for the considered class of problems \cite{Hazan-Kale}.
\end{remark}

\begin{remark}\label{Rem_many_constr}
	Let's consider a modification of the proposed Algorithm 2 for the case of a set of functional constraints $g_m: Q\rightarrow \mathbb{R}$ ($m = \overline{1, K}$). We assume, that all the functionals $g_m$ satisfy the Lipschitz condition:
	\begin{equation}\label{eq1}
	|g_m(x)-g_m(y)|\leq M||x-y||\;\forall x,y\in Q, \; m = \overline{1, K}.
	\end{equation}
	In this case, instead of a set of convex functional constraints $\{g_m(\cdot)\}_{m=1}^{K}$ we can consider one constraint, given as $g: Q \rightarrow \mathbb{R}$, where
	$$
	g(x) = \max\limits_{m = \overline{1, K}} g_m(x), \quad |g(x)-g(y)|\leq M||x-y||\;\forall x,y\in Q.
	$$
	This method will be also optimal, but in practice it can give better accuracy (see Remark \ref{Ston_Rem_4} below).
\end{remark}

\begin{algorithm}[t]
	\begin{algorithmic}[1]
		\caption{Online Optimization: Adaptive Mirror Descent Algorithm Modification for the Case of Many Constraints}
		\label{alg03}
		\REQUIRE  $\varepsilon, N, \Theta_0^2,Q,d(\cdot), x^0$
		\STATE $i:=1, \; k:=0;$
		\REPEAT
		\IF{$g(x^k) \leqslant \varepsilon$}
		\STATE $M_k := \|\nabla f_i(x^k)\|_*;$
		\STATE $h_k = \Theta_0 \left(\sum\limits_{t = 0}^k M_t^2\right)^{-1/2};$
		\STATE $x^{k+1}:=Mirr[x^k](h_k\nabla f_i(x^k));$
		\STATE $i:=i+1;$
		\STATE $k:=k+1;$
		\ELSE
		\STATE $M_k := \|\nabla g_{m(k)}(x^k)\|_*$ for some $g_{m(k)}(\cdot))$: $g_{m(k)}(x^k))> \varepsilon$
		\STATE $h_k = \Theta_0 \left(\sum\limits_{t = 0}^k M_t^2\right)^{-1/2};$
		\STATE $x^{k+1}:=Mirr[x^k](h_k\nabla g_{m(k)}(x^k));$
		\STATE $k:=k+1;$
		\ENDIF
		\UNTIL{$i=N+1$}
		\STATE Guaranteed accuracy:
		\begin{equation} \label{eqiv_crit_2}
		\delta:= \frac{2\Theta_0}{N} \left(\sum\limits_{i=0}^{N+N_J-1} M_i^2\right)^{1/2} - \varepsilon \cdot \frac{N_J}{N}.
		\end{equation}
	\end{algorithmic}
\end{algorithm}

\section{The Case of Negative Regret}

Now we consider the situation, when after the stopping of any of the above algorithms, it turns out that the regret is negative.
In this case the following inequality
\begin{equation}\label{equiv_negat_regret}
\frac{1}{N}\sum\limits_{i=1}^Nf_i(x^k)-\min\limits_{x\in Q}\frac{1}{N}\sum\limits_{i=1}^Nf_i(x) \leq 0
\end{equation}
holds. It is already impossible to justify the optimality of the number of non-productive steps in view of the right-hand side of inequality \eqref{equiv_negat_regret}.

Note that the set of productive steps is not empty, because for arbitrary $p$ steps when the inequality
$$\sum\limits_{k=1}^p\frac{1}{M_k^2}\geq \frac{2\Theta_0^2}{\varepsilon^2}$$
is satisfied, one of these $p$ steps will necessarily be productive (see \cite{bib_Adaptive,bib_Stonyakin}). If all the other $p-1$ steps are non-productive (without loss of generality let the last step be productive), then
$$\sum\limits_{k=1}^{p-1}\frac{1}{M_k^2}< \frac{2\Theta_0^2}{\varepsilon^2}$$
and
$$p-1<\frac{2M^2\Theta_0^2}{\varepsilon^2}.$$

It is clear, that running the method for a sufficiently long time, it is possible to achieve $N$ productive steps. At the same time between each two successive productive steps there will be no more than $\frac{2M^2\Theta_0^2}{\varepsilon^2}$ non-productive steps, i.e. the number of all non-productive steps will be no more than
$$\frac{2M^2\Theta_0^2}{\varepsilon^2}N.$$

In comparison with the previous items, for $ \varepsilon = \frac {C} {\sqrt {N}} $ there will be no more than
\begin{equation}\label{Non_neg_nonprod}
\frac{2M^2\Theta_0^2}{\varepsilon^2}N=O(N^2)
\end{equation}non-productive steps.

\section{Numerical Experiments}
To compare of Algorithms \ref{alg01}, \ref{alg02} and \ref{alg03}, some numerical tests were carried out. Consider four different examples with objective function
$$f(x) = \frac{1}{N} \displaystyle \sum_{i=1}^{N} \left| \langle a_i , x \rangle - b_i \right|. $$
For the coefficients $a_i$ and constants $b_i$ for $i=1,\dots , N$, with different values of $N$. Let $A \in \mathbb{R}^{N \times 11} $ be a matrix with entries drawn from different random distiributions. Then $a_i^T$ are rows in the matrix $A' \in \mathbb{R}^{N \times 10} $, which is introduced from $A$, by eliminating the last column, and $b_i$ are the entries of the last column in the matrix $A$. In details, entries of $A$ drawn
\begin{itemize}
	\item In example 1, from a normal distribution with mean (center) equalling $0$ and standard deviation (width) equalling $1$.
	\item In example 2, from a uniform distribution over $[0,1)$.
	\item In example 3, from the standard exponential distribution with a scale parameter of $1$.
	\item In example 4, from a Gumbel distribution with the location of the mode equalling $1$ and the scale parameter equalling $2$.
\end{itemize}
For the function of constraints $g(x) = \max\limits_{i \in \overline{1, m}} g_i(x)$, we take $m = 3$ and the functionals $g_i(x) =\langle \alpha_i , x \rangle $, where $\alpha_i^{T}$ are the rows of the matrix \\
$$
\left(
\begin{array}{cccccccccc}
1 & 1  & 1  & 1  & 1  & 1  & 1  & 1  & 1  & 1\\
1\, & 2\, & 3\, & 4\, & 5\, & 6\, & 7\, & 8\, & 9\, & 10\\
1 & 2 & 4 & 6 & 8 & 10 & 12 & 14 & 16 & 18
\end{array}
\right)
$$
We choose standard Euclidean proximal setup as a prox-function, starting point $x^{0} = \dfrac{(1,1,...,1)}{\sqrt{10}}$, $\Theta_0 = 3$, $\varepsilon = \frac{1}{\sqrt{N}}$ and
$$Q = \{x = (x_1, x_2, ..., x_{10}) \in \mathbb{R}^{10}\,|\, x_1^2 + x_2^2 + ... + x_{10}^2 \leq 1\}.$$
The results of the work of Algorithms \ref{alg01}, \ref{alg02} and \ref{alg03} are represented in Table \ref{tab1}, Table \ref{tab2} and Table \ref{tab3} below, respectively, demonstrate the comparison between these algorithms. The number of non-productive steps are denoted by {\it nonprod.}, time is given in seconds and parts of the second, $\delta$ is guaranteed accuracy of the solution approximation found (sequence $\{x^k\}_{k\in I}$ on productive steps).

All experiments were implemented in Python 3.4, on computer fitted with Intel(R) Core(TM) i7-8550U CPU @ 1.80GHz, 1992 Mhz, 4 Core(s), 8 Logical Processor(s). RAM of the computer is 8GB.

\begin{table}[]
	\centering
	\caption{Results of Algorithm \ref{alg01}.}
	\label{tab1}
	\begin{tabular}{|c|c|c|c|}
		\hline
		\multirow{3}{*}{} & nonprod.    & time & $\delta$     \\ \cline{1-4}
		ex. 1, $N=3000$   &7041  &00.444 &187.473    \\ \hline
		ex. 2, $N=6000$   &12645 &00.812 &132.565   \\ \hline
		ex. 3, $N=7000$   &15814 &00.958 &122.730   \\ \hline
		ex. 4, $N=10000$  &24971 &01.523 &102.682   \\ \hline
	\end{tabular}
\end{table}

\begin{table}[]
	\centering
	\caption{Results of Algorithm \ref{alg02}.}
	\label{tab2}
	\begin{tabular}{|c|c|c|c|}
		\hline
		\multirow{3}{*}{} & nonprod.    & time & $\delta$     \\ \cline{1-4}
		ex. 1, $N=3000$   &39    &00.149 &0.426    \\ \hline
		ex. 2, $N=6000$   &2821  &00.404 &0.223   \\ \hline
		ex. 3, $N=7000$   &5543  &00.586 &0.405   \\ \hline
		ex. 4, $N=10000$  &12576 &01.104 &0.692   \\ \hline
	\end{tabular}
\end{table}

From Table \ref{tab1} and Table \ref{tab2} one can see, that the adaptive Algorithm \ref{alg02} always works better than non-adaptive Algorithm \ref{alg01}. It is clearly shown in all the examples by the number of non-productive steps, running time of the algorithms and guaranteed accuracy $\delta$. Where the number of non-productive steps and $\delta$ produced by Algorithm \ref{alg02} is very small compared to the Algorithm \ref{alg01}.

From Table \ref{tab3}, we can see, that there is a difference between the number of non-productive steps produced by Algorithms \ref{alg02} and \ref{alg03}, but the guaranteed accuracy $\delta$ and the running time produced by Algorithm \ref{alg03} is smaller compared to Algorithm \ref{alg02}.

\begin{table}[]
	\centering
	\caption{Results of Algorithm \ref{alg03}.}
	\label{tab3}
	\begin{tabular}{|c|c|c|c|}
		\hline
		\multirow{3}{*}{} & nonprod.    & time & $\delta$     \\ \cline{1-4}
		ex. 1, $N=3000$   &47   &00.121 &0.414   \\ \hline
		ex. 2, $N=6000$   &2835 &00.333 &0.220   \\ \hline
		ex. 3, $N=7000$   &5563 &00.454 &0.394   \\ \hline
		ex. 4, $N=10000$  &12885&00.807 &0.680   \\ \hline
	\end{tabular}
\end{table}

\begin{remark}\label{Ston_Rem_4}
To show the advantages of Algorithm \ref{alg03}, as compared to Algorithm \ref{alg02}, one additional numerical test was carried out.
Let's now take the functionals of constraints $g_i$, $i=1,2,3$ as follows
$$
	g_1(x) = \displaystyle \sum_{i=1}^{10}i \cdot x_i +1 , \ \  g_2(x) = \displaystyle \sum_{i=1}^{10}10i \cdot x_i , \ \ g_3(x) = \displaystyle \sum_{i=1}^{10}50i \cdot x_i .
$$
with the same all previous parameters: starting point $x^{0} = \dfrac{(1,1,...,1)}{\sqrt{10}}$, $\Theta_0 = 3$,
$$Q = \{x = (x_1, x_2, ..., x_{10}) \in \mathbb{R}^{10}\,|\, x_1^2 + x_2^2 + ... + x_{10}^2 \leq 1\},$$
but with $\varepsilon = 0.5$. Table \ref{tab4} below demonstrate the comparison between Algorithms \ref{alg02} and \ref{alg03}, for the objective function $f(x)= \frac{1}{3} \displaystyle \sum_{i=1}^{3} f_{i}(x)$, where
$$
	f_1(x) = \sqrt{\displaystyle \sum_{i=1}^{9}( x_i + x_{i+1} )^2}, \ f_2(x) =  \sqrt{0.1\left(  \displaystyle\sum_{i=1}^{10} x_i^2 + \displaystyle\sum_{i=1}^{9} x_i x_{i+1} \right)}, \ f_3(x) = \sqrt{\displaystyle\sum_{i=1}^{10} x_i^2}.
$$
	
	\begin{table}[]
		\centering
		\caption{Results of algorithms \ref{alg02} and \ref{alg03}.}
		\label{tab4}
		\begin{tabular}{|c|c|c|c|}
			\hline
			\multirow{3}{*}{} ex. 5, $N=3$ & nonprod.    & time & $\delta$     \\ \cline{1-4}
			Algorithm 2   &1 &00.044 &1961.954    \\ \hline
			Algorithm 3   &2 &00.030 &9.608  \\ \hline
		\end{tabular}
	\end{table}
	
	From Table \ref{tab4}, one can see, that Algorithm \ref{alg03} works better than Algorithm \ref{alg02}, since the difference between the non-productive steps is very small, equalling only one, and the guaranteed accuracy $\delta$ produced by Algorithm \ref{alg03} is very small compared to the precision produced by Algorithm \ref{alg02}.
	
\end{remark}


\begin{thebibliography}{30}
	%
	
	\bibitem{bayandina2018primal-dual}
	Bayandina, A., Gasnikov, A., Gasnikova, E., Matsievsky, S.:
	Primal-dual mirror descent for the stochastic programming problems with functional constraints.
	Computational Mathematics and Mathematical Physics. (accepted) (2018). \url{https://arxiv.org/pdf/1604.08194.pdf}. (in Russian)
	
	\bibitem{bib_Adaptive}
	Bayandina, A., Dvurechensky, P., Gasnikov, A., Stonyakin, F., Titov, A.:
	Mirror descent and convex optimization problems with non-smooth inequality constraints.
	In: LCCC Focus Period on Large-Scale and Distributed Optimization. Sweden, Lund: Springer. (accepted) (2017).
	\url{https://arxiv.org/abs/1710.06612}
	
	\bibitem{beck2010comirror}
	Beck, A., Ben-Tal, A., Guttmann-Beck, N., Tetruashvili, L.:
	The comirror algorithm for solving nonsmooth constrained convex problems. Operations Research Letters 38(6), pp. 493--498 (2010).
	
	\bibitem{beck2003mirror}
	Beck, A., Teboulle, M.: Mirror descent and nonlinear projected subgradient methods for convex optimization.
	Oper. Res. Lett. 31(3), pp. 167--175 (2003).
	
	\bibitem{bib_Nemirovski}
	Ben-Tal, A., Nemirovski, A.: Lectures on Modern Convex Optimization.
	Society for Industrial and Applied Mathematics, Philadelphia (2001).
	
	\bibitem{bib_ttd}
	Ben-Tal, A., Nemirovski, A.: Robust Truss Topology Design via semidefinite programming.
	SIAM Journal on Optimization vol. 7, no. 4, pp. 991--1016 (1997).
	
	\bibitem{bib_Boyd}
	Boyd, S., Vandenberghe, L.: Convex Optimization. Cambridge University Press, New York (2004).
	
	\bibitem{Bubeck-Eldan}
	Bubeck S., Eldan R.: Multi-scale exploration of convex functions and bandit convex optimization. e-print, (2015). \url{http://research.microsoft.com/en-us/um/people/sebubeck/ConvexBandits.pdf}
	
	\bibitem{Bubeck-Bianchi}
	Bubeck S., Cesa-Bianchi N.: Regret analysis of stochastic and nonstochastic multi-armed bandit problems. Foundation and Trends in Machine Learning, vol. 5, no. 1, pp. 1--122 (2012).
	
		
\bibitem{Gasnikov2015}
	Gasnikov A.V., Lagunovskaya A.A.,  Morozova L.E.: On the relationship between simulation logit dynamics in the population game theory and a mirror descent method in the online optimization using the example of the shortest path problem. PROCEEDINGS OF MIPT, vol. 7, no. 4, pp. 104--113 (2015). (in Russian)
	
\bibitem{Gasnikov2017}
	Gasnikov A.V., Lagunovskaya A.A., Usmanova I.N., Fedorenko F.A., Krymova E.A.: Stochastic online optimization. Single-point and multi-point non-linear multi-armed bandits. Convex and strongly-convex case. Automation and Remote Control, vol. 78, Issue 2, pp. 224--234, (2017).

	\bibitem{Hazan-Kale}
	Hazan E., Kale S.: Beyond the regret minimization barrier: Optimal algorithms for stochastic strongly-convex optimization. JMLR. vol. 15, pp. 2489--2512 (2014).
	
	\bibitem{Hazan-2016}
	Hazan E. Introduction to online convex optimization. Foundations and Trends in Optimization, vol. 2, no. 3--4, pp. 157--325 (2015).
	
\bibitem{Jenatton_adaptive_with_constraints} Jenatton R., Huang J., Archambeau C.: Adaptive Algorithms for Online Convex Optimization
		with Long-term Constraints. (2015). \url{https://arxiv.org/abs/1512.07422}.
	
\bibitem{Kalai_efficient_algorithms}
Kalai A., Vempala S.: Efficient algorithms for online decision problems. Journal of Computer and System Sciences, vol. 71, pp. 291 -- 307 (2005).

	\bibitem{Lugosi-Bianchi}
	Lugosi G., Cesa-Bianchi N.: Prediction, learning and games. New York, Cambridge University Press, (2006).
	
	\bibitem{nemirovskii1979efficient}
	Nemirovskii, A.:
	Efficient methods for large-scale convex optimization problems. Ekonomika i Matematicheskie Metody (1979). (in Russian)
	
	\bibitem{nemirovsky1983problem}
	Nemirovsky, A., Yudin, D.: Problem Complexity and Method Efficiency in Optimization. J. Wiley \& Sons, New York (1983).
	
	\bibitem{bib_Nesterov}
	Nesterov, Y.: Introductory Lectures on Convex Optimization: A Basic Course. Kluwer Academic Publishers, Massachusetts (2004).
	
	\bibitem{polyak1967general}
	Polyak, B.: A general method of solving extremum problems. Soviet Mathematics Doklady, vol. 8, no. 3, pp. 593--597 (1967). (in Russian)
	
	\bibitem{shor1967generalized}
	Shor, N. Z.: Generalized gradient descent with application to block programming. Kibernetika vol. 3, no. 3, pp. 53--55 (1967). (in Russian)
	
	\bibitem{bib_Stonyakin}
	F.S. Stonyakin, M. S. Alkousa, A. N. Stepanov, M. A. Barinov.: Adaptive mirror descent algorithms
	in convex programming problems with Lipschitz constraints. Trudy Instituta Matematiki i Mekhaniki URO RAN, vol. 24, no. 2, pp. 266 -- 279 (2018).
	
	\bibitem{bib_Shpirko}
	Shpirko, S., Nesterov Y.: Primal-dual subgradient methods for huge-scale linear conic problem.
	SIAM Journal on Optimization, vol. 24, no. 3, pp. 1444 -- 1457 (2014).
	
	\bibitem{bib_Vasilyev}
	Vasilyev, F.: Optimization Methods. Fizmatlit, Moscow (2002). (in Russian)
	
	
	
\end{thebibliography}
\end{document}